\documentclass{amsart}

\usepackage{graphicx}
\usepackage{amsmath, amssymb}
\usepackage{enumerate}
\newcommand{\abs}[1]{\lvert#1\rvert}

\begin{document}

\title[Discrete differential geometry of proteins]
{Discrete differential geometry of proteins: 
a new method for encoding three-dimensional structures of proteins}

%    Information for the author
\author{Naoto Morikawa}
\address{GENOCRIPT, 27-22-1015, Sagami-ga-oka 1-chome, Zama-shi, Kanagawa 228-0001 Japan}
\email{nmorika@f3.dion.ne.jp}

%    General info
\date{Posted June 5, 2005.}

\begin{abstract}
In nature the three-dimensional structure of a protein is encoded in the corresponding gene. 
In this paper we describe a new method for encoding the three-dimensional structure of a 
protein into a binary sequence. The feature of the method is the correspondence between 
protein-folding and ``integration''. A protein is approximated by a folded tetrahedron sequence. 
And the binary code of a protein is obtained as the ``second derivative'' of the shape of the 
folded tetrahedron sequence. With this method at hand, we can extract static structural 
information of a protein from its gene. And we can describe the distribution of three-dimensional 
structures of proteins without any subjective hierarchical classification. 
\end{abstract}

\maketitle

\setcounter{tocdepth}{3}
\tableofcontents

%----------------------------------------------------------------------
\section{Overview}

\begin{figure}[b]
\includegraphics{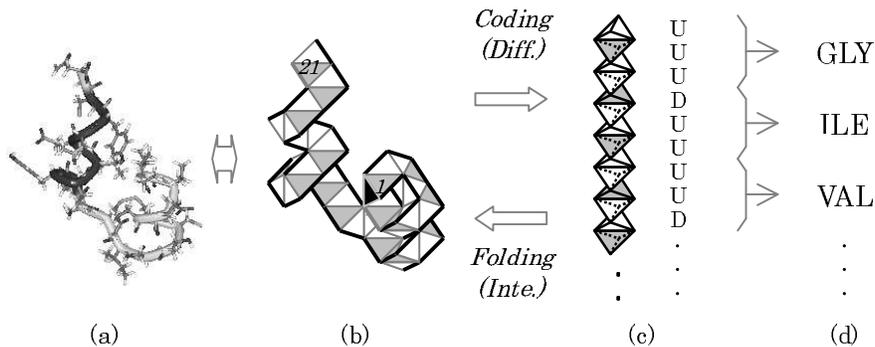}
\caption{Overview. (a): Schematic diagram of 2HIU chain A (Insulin, human). 
(b): Approximation of 2HIU by a tetrahedron sequence. 
(c): The $U/D$ sequence of approximation (b). 
(d): The amino-acid sequence of 2HIU. 
(The figure (a) is prepared using \textbf{WebLab Viewer} (Molecular Simulations Inc.).)}
\label{figHIU}
\end{figure}

In nature  the three-dimensional structure of a protein is encoded in the corresponding gene. 
In this paper we describe a new method for encoding the three-dimensional structure of a 
protein into a binary sequence (Fig.\ref{figHIU}).

In the method a protein is approximated by a tetrahedron sequence. For example, 
approximation Fig.\ref{figHIU}(b) is obtained by folding tetrahedron sequence 
Fig.\ref{figHIU}(c), where three tetrahedrons are assigned for each amino-acid. 
We would obtain more precise approximation if we use more tetrahedrons.

The feature of the method is the correspondence between protein-folding and ``integration''. 
And the binary sequence is obtained as the ``second derivative'' of the shape of the folded 
tetrahedron sequence.

With this method at hand, we can extract static structural information of a protein from its gene. 
And we can describe the distribution of three-dimensional structures of proteins 
without any subjective hierarchical classification.

\section{Basic idea: encoding of two-dimensional objects}

\begin{figure}[tb]
\includegraphics{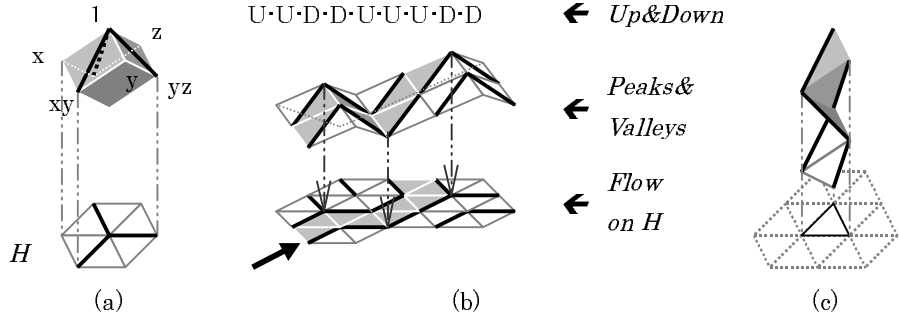}
\caption{Basic idea. 
(a): Unit cube in $\mathbb{R}^3$ and its projection on $H$.
(b): Slant-tile sequences and flat-tile sequences defined by $Cone^\ast \{x^2/z, 1 \}$. 
(c): Slant-tiles over a flat-tile on $H$.}
\label{fig1}
\end{figure}

For simplicity we shall explain the basic idea behind the paper in the case of 
two-dimensional objects, where we use triangle sequences for approximation.

\subsection{Triangle sequence}

Consider a unit cube in the three-dimensional Euclidean space $\mathbb{R}^3$ whose vertices 
are given by $v_{1}$, $v_{x}$, $v_{y}$, $v_{z}$, $\ldots$, and $v_{xyz}$, 
where $v_{x^ly^mz^n}:=(l,m,n)\in \mathbb{Z}^3$ (Fig.\ref{fig1}(a)). And draw lines 
$\overline{v_{1}v_{xy}}$, $\overline{v_{1}v_{yz}}$ and $\overline{v_{1}v_{xz}}$. 
Then, each of three upper faces is divided into two slant-triangle-tiles. 
For example, $v_{1}v_{x}v_{xy}v_{y}$ is divided into two slant-tiles $v_{1}v_{x}v_{xy}$ and 
$v_{1}v_{y}v_{xy}$.

Firstly, by piling up these unit cubes in the direction from $v_{xyz}$ to $v_{1}$, we obtain 
``peaks and valleys'' with a ``drawing'' on it. The drawing is uniquely determined by its 
peaks and divides the surface into a collection of slant-triangle-tile sequences. For example, 
the drawing of Fig.\ref{fig1}(b) is determined by two peaks, left $(2,0,-1)$ and right 
$(0,0,0)$. And we denote the drawing by $Cone^\ast \{x^2/z, 1 \}$.

Secondly, by the projection onto the hypersurface $H:=\{(a,b,c) \in \mathbb{R}^3 \ |\ a+b+c=0\}$, 
we obtain a division of $H$ into a collection of flat-triangle-tile sequences.  
For example, the gray slant-tile sequence is projected onto the gray flat-tile sequence on $H$ 
in Fig.\ref{fig1}(b). We write $a[uv]$ for slant-tile $v_{a}v_{au}v_{auv}$ and $\abs{a[uv]}$ 
for the corresponding flat-tile. For example, $1[xy]$ for $v_{1}v_{x}v_{xy}$. 
Note that there are three types of slant-tiles over a flat-tile (Fig.\ref{fig1}(c)). 
We shall see in the appendix that ``peaks and valleys'' specifies a ``discrete vector field'' 
of flat-tiles on $H$.

Finally we obtain a binary code of the shape of a flat-tile sequence by arranging up ($U$) and down 
($D$) of the corresponding slant-tile sequence. For example, the gray flat-tile sequence in 
Fig.\ref{fig1}(b) is encoded into $U/D$ sequence
\[
U-U-D-D-U-U-U-D-D. 
\]
In general we need more than one drawing to encode a flat-tile sequence because of overlaps 
among its peaks (Fig.\ref{fig2}(c)). Each drawing encodes a part of the flat-tile sequence and its 
code is obtained by patching those ``local codes'' together. 

\subsection{Encoding of two-dimensional objects}

\begin{figure}[tb]
\includegraphics{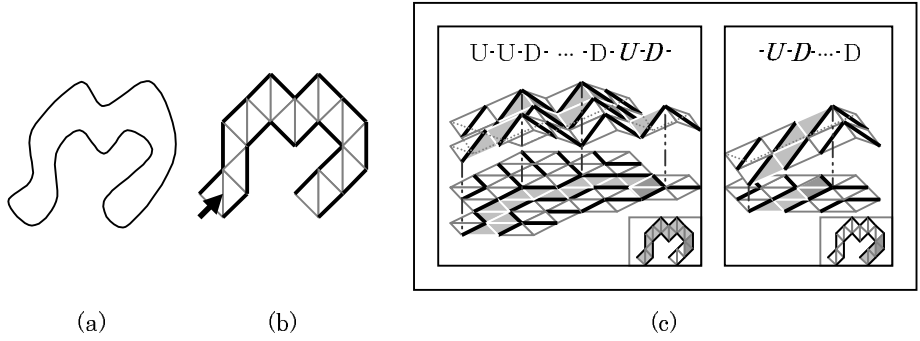}
\caption{Encoding of a two-dimensional object. (a): Two-dimensional object. 
(b): Approximation by a triangle sequence. 
(c): Two drawings $Cone^\ast \{1, z/y^2, z^2/(xy^2), z^3/x^2 \}$ and 
$Cone^\ast \{ z^3/x^2 \}$ which encode approximation (b). }
\label{fig2}
\end{figure}

Now let's encode the two-dimensional object shown in Fig.\ref{fig2}(a). 
First of all we should give a flat-tile sequence which approximates the object 
(Fig.\ref{fig2}(b)). 

Then, using encoding table Table \ref{tab1}(a), we obtain a binary code of the object 
(Fig.\ref{fig2}(c)). The process is going on as follows: 
\begin{enumerate}[\ \ \ \ \ \ \ \ Step 1.]
\item   Choose an initial value, say $U$, 
\item   By the second row of the table, the second value is $U$, 
\item   By the fourth row of the table, the third value is $D$, \ldots.
\end{enumerate}
As the result we obtain $U/D$ sequence 
\begin{equation}
U-U-D-D-U-U-U-D-D-U-U-D-D-D-U-D-D-D-D. \label{seq: basic}
\end{equation}

\begin{table}[tb]
\includegraphics[scale=0.7]{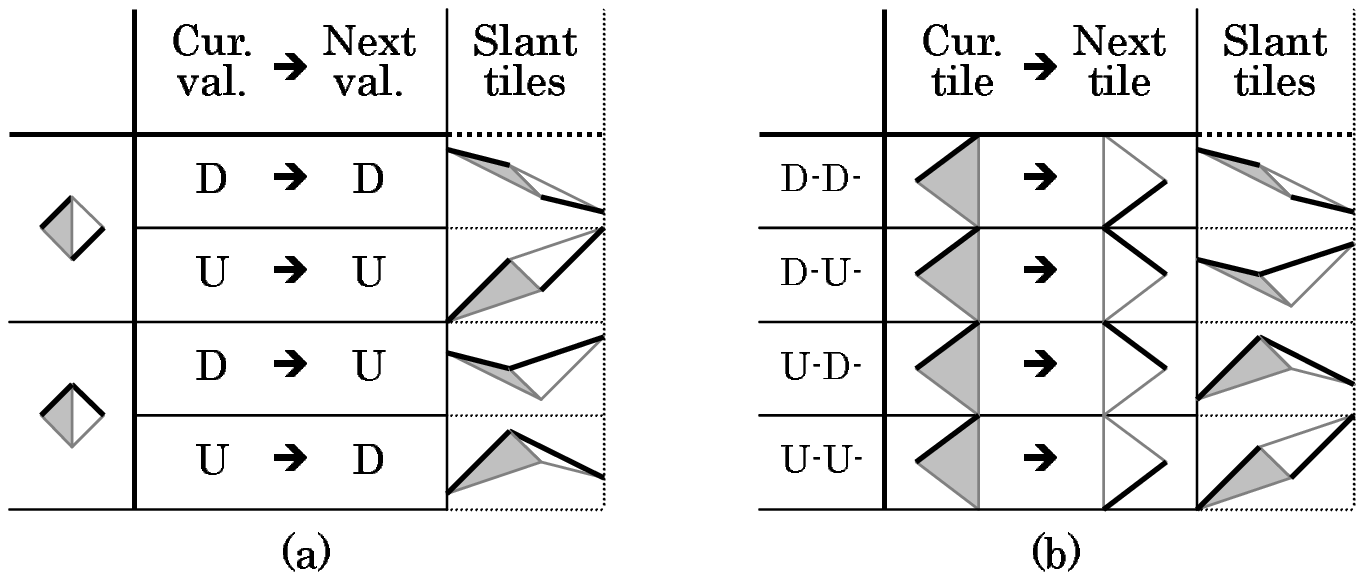}
\caption{Tables for two-dimensional objects. (a): Encoding table. 
(b): Decoding table. (The gray tile is the current one.)}
\label{tab1}
\end{table}

Fig.\ref{fig2}(c) shows the corresponding slant-tile sequence. In this case we need two 
drawings because of the overlap between two peaks $z^2/(y^2x)$ and $z^3/x^2$. 
The left drawing $Cone^\ast \{1, z/y^2, z^2/(xy^2), z^3/x^2 \}$ corresponds to the first sixteen tiles 
and the right drawing $Cone^\ast \{ z^3/x^2 \}$ to the last five tiles.

\subsection{Decoding of $U/D$ sequences in $\mathbb{R}^2$}

To decode $U/D$ sequences in $\mathbb{R}^2$ we use decoding table Table \ref{tab1}(b). 
For example, decoding process of $U/D$ sequence (\ref{seq: basic}) is going on as follows: 
\begin{enumerate}[\ \ \ \  \ \ \ \ Step 1.]
\item   Choose an initial flat-tile, say $\abs{x[yx]}$, 
\item   By the fourth row of the table, the second flat-tile is $\abs{1[xy]}$,
\item   By the third row of the table, the third flat-tile is $\abs{1[xz]}$, \ldots.
\end{enumerate}
As the result we obtain the flat-tile sequence shown in Fig.\ref{fig2}(b).

\section{Encoding of three-dimensional objects}

If we consider unit cubes in the four-dimensional Euclidean space $\mathbb{R}^4$, 
we shall obtain a three-dimensional drawing made up of slant-``tetrahedron''-tiles. 
And we approximate a three-dimensional object by a tetrahedron sequence 
(Fig.\ref{figHIU}(c)), where
\begin{enumerate}
\renewcommand{\labelenumi}{(\arabic{enumi})}
\item each tetrahedron consists of four short edges and two long edges, 
       where the ratio of the length is $\sqrt{3}/2$  and
\item successive tetrahedrons are connected via a long edge and have the rotational 
       freedom around the edge.
\end{enumerate}

\subsection{Tetrahedron sequence}

\begin{figure}[tb]
\includegraphics{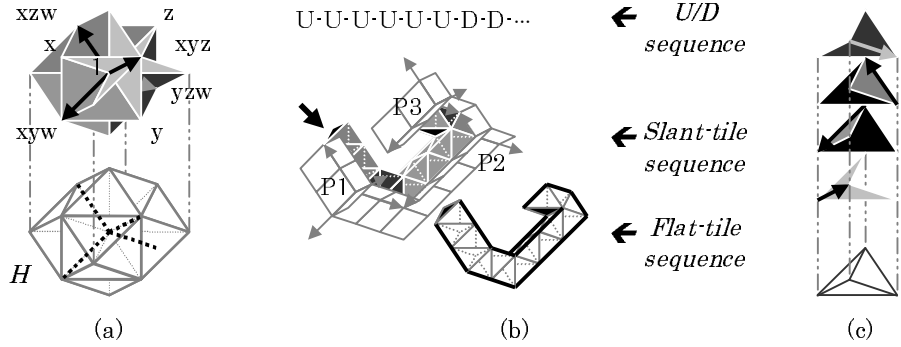}
\caption{Encoding of three-dimensional objects. 
(a): Unit cube in $\mathbb{R}^4$ and its projection on $H$.
(b): Slant-tile sequence and flat-tile sequence defined by three peaks $P_1=1$, 
$P_2=z^2/(x^3yw)$, and $P_3=z^2/(x^2y^2w)$.  
(c): Slant-tiles over a flat-tile on $H$.
(In the figures arrows indicate the direction of ``down''.) }
\label{fig3}
\end{figure}

Consider a unit cube in the four-dimensional Euclidean space $\mathbb{R}^4$ whose vertices 
are given by $v_{1}$, $v_{x}$, $v_{y}$, $v_{z}$, $v_{w}$, $\ldots$, and $v_{xyzw}$, 
where $v_{x^ly^mz^nw^k}:=(l,m,n,k)\in \mathbb{Z}^4$ (Fig.\ref{fig3}(a)). 
And divide each of four upper three-dimensional faces into six slant-tetrahedron-tiles. 
For example, the face defined by $v_{1}$, $v_{x}$, $v_{z}$, and $v_{y}$ is divided into 
six slant-tiles $v_{1}v_{x}v_{xy}v_{xyz}$, $v_{1}v_{y}v_{yx}v_{yxz}$, $v_{1}v_{y}v_{yz}v_{yzx}$, 
$v_{1}v_{z}v_{zy}v_{zyx}$, $v_{1}v_{x}v_{xz}v_{xzy}$, and $v_{1}v_{z}v_{zx}v_{zxy}$.

Firstly, by piling up these unit cubes in the direction from $v_{xyzw}$ to $v_{1}$, we obtain 
four-dimensional ``peaks and valleys'' with a three-dimensional ``drawing'' on it. 
The drawing is uniquely determined by its peaks and divides the three-dimensional surface into 
a collection of slant-tetrahedron-tile sequences. For example, 
the drawing of Fig.\ref{fig3}(b) is determined by three peaks $P_1$, $P_2$ and $P_3$. 
And we denote the drawing by $Cone^\ast \{P_1, P_2, P_3 \}$.

Secondly, by the projection onto the hypersurface 
$H:=\{(a,b,c,d) \in \mathbb{R}^4 \ |\  a+b+c+d=0\}$, we obtain a 
division of $H$ into a collection of flat-tetrahedron-tile sequences.  
For example, Fig.\ref{fig3}(b) shows a slant-tile sequence and its projection onto $H$. 
We write $a[uvw]$ for slant-tile $v_{a}v_{au}v_{auv}v_{auvw}$ and $\abs{a[uvw]}$ 
for the corresponding flat-tile. 
Note that there are four types of slant-tiles over a flat-tile (Fig.\ref{fig3}(c)). 
For example, $\abs{1[xyz]}=\abs{x[yzw]}=\abs{xy[zwx]}=\abs{xyz[wxy]}$. 

Finally we obtain a binary code of the shape of a flat-tile sequence by arranging up ($U$) and down 
($D$) of the corresponding slant-tile sequence. For example, the flat-tile sequence shown in 
Fig.\ref{fig3}(b) is encoded into $U/D$ sequence
\begin{equation}
U^6-D^4-U^7-D-U-D.  \label{seq: three-dim}
\end{equation}

\subsection{Encoding of three-dimensional objects}

To encode three-dimensional objects we use encoding table Table \ref{tab2}(a). 
For example, encoding of the flat-tile sequence shown in Fig.\ref{fig3}(b) proceeds as follows: 
\begin{enumerate}[\ \ \ \ \ \ \ \ Step 1.]
\item   Choose an initial value, say $U$, 
\item   By the second row of the table, the second value is $U$, 
\item   By the second row of the table, the third value is $U$, \ldots.
\end{enumerate}
As the result we obtain $U/D$ sequence (\ref{seq: three-dim}).

\subsection{Decoding of $U/D$ sequences in $\mathbb{R}^3$}

To decode $U/D$ sequences in $\mathbb{R}^3$ we use decoding table Table \ref{tab2}(b). 
For example, decoding of $U/D$ sequence (\ref{seq: three-dim}) proceeds as follows: 
\begin{enumerate}[\ \ \ \  \ \ \ \ Step 1.]
\item   Choose an initial flat-tile, say $\abs{xw^2z^2[xzw]}$, 
\item   By the fourth row of the table, the second flat-tile is $\abs{xwz^2[wxz]}$,
\item   By the fourth row of the table, the third flat-tile is $\abs{xwz[zwx]}$, \ldots.
\end{enumerate}
As the result we obtain the flat-tile sequence shown in Fig.\ref{fig3}(b).

\begin{table}[tb]
\includegraphics[scale=0.7]{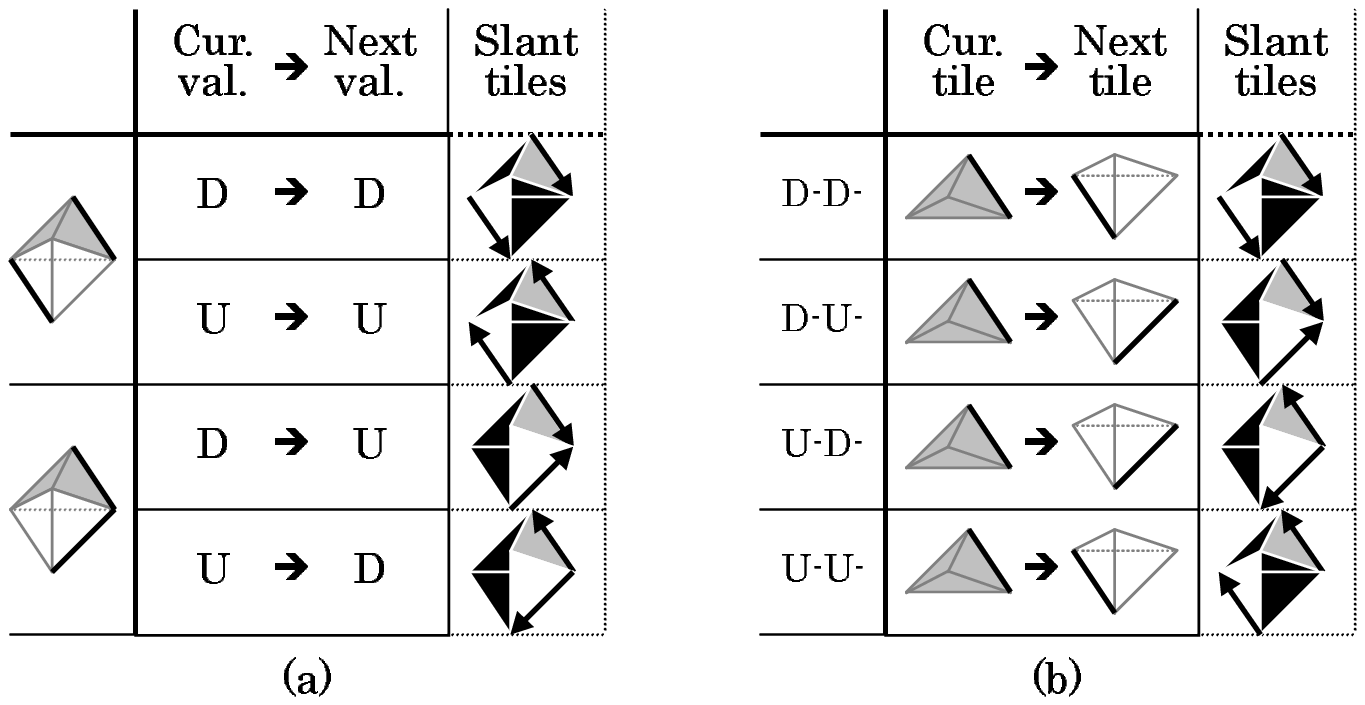}
\caption{Tables for three-dimensional objects. 
(a): Encoding table. (b): Decoding table. (The gray tile is the current one.)}
\label{tab2}
\end{table}

\section{Examples}

\subsection{Double helix}

\begin{figure}[tb]
\includegraphics{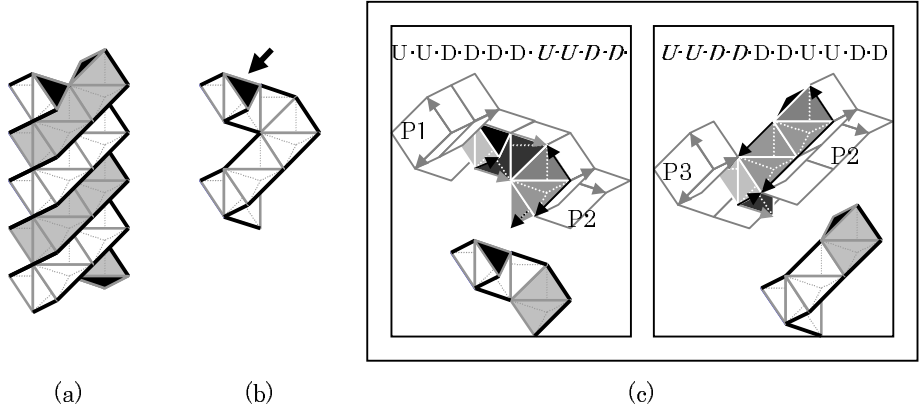}
\caption{Double helix. 
(a): Double helix formed by two tetrahedron sequences. 
(b): Part of the helix. (c): Slant-tile sequence and flat-tile sequence defined 
by $P_1=1$, $P_2=y^2z/x$, and $P_3=y^2w^2$.}
\label{figHelix}
\end{figure}

Here let's consider the double helix shown in Fig.\ref{figHelix}(a) which has $12$ tiles per turn. 
(Cf. DNA has an average of $10.9$ (type A) or  $10$ (type B) nucleotide pairs per turn (\cite{BT}).) 
To encode the shape of the helix, it is enough to consider the flat-tile sequence shown 
in Fig.\ref{figHelix}(b).

Using Table \ref{tab2}(a) with initial slant-tile $y[zxy]$, we obtain two drawings of 
Fig.\ref{figHelix}(c). $Cone^\ast \{P_1, P_2\}$ (left) encodes the first ten tiles. 
And $Cone^\ast \{P_2, P_3\}$ (right) encodes the last ten tiles. By patching these 
local codes together, we obtain the $U/D$ code of helix Fig.\ref{figHelix}(b): 
\[
U-U-D-D-D-D\ - \ U-U-D-D\ - \ D-D-U-U-D-D.
\]

\subsection{2HIU chain A (Insulin, human)}

Next let's consider the three-dimensional structure of 2HIU chain A (Fig.\ref{figHIU}). 
Using Table \ref{tab2}(a) with initial slant-tile $zw[xyz]$, we obtain eight drawings: 
\begin{align*}
&\text{-\ \ } Cone^\ast \{ z/y, 1/(x^2w), 1/(x^2z)\} \quad &\text{for $[1,14]$,\ } \\
&\text{-\ \ } Cone^\ast \{ 1/(xy), 1/(x^2w), 1/(x^2z)\} \quad &\text{for $[7,18]$,\ } \\
&\text{-\ \ } Cone^\ast \{ 1/(xy), 1/(x^3zw), 1/(x^3z^2), w/(xyz)\} \quad &\text{for $[13,29]$,} \\
&\text{-\ \ } Cone^\ast \{ 1/(xyz^2), 1/(x^3zw), 1/(x^3z^2), xw/y^2\} \quad &\text{for $[16,42]$,} \\
&\text{-\ \ } Cone^\ast \{ xw^2/(yz), w/y, xw/y^2\} \quad &\text{for $[36,45]$,} \\
&\text{-\ \ } Cone^\ast \{ xw^2/(yz), 1/y^2, x/y^3\} \quad &\text{for $[40,51]$,} \\
&\text{-\ \ } Cone^\ast \{ x/(y^4z), 1/y^2, x/(y^4w)\} \quad &\text{for $[45,57]$,} \\
&\text{-\ \ } Cone^\ast \{ x/(y^4z), 1/(y^4w^2)\} \quad &\text{for $[52,63]$.}
\end{align*}
($[n,m]$ denotes the part of the sequence from the $n$-th tile to the $m$-th tile.)

By patching these local codes together, we obtain the $U/D$ code of the three-dimensional 
structure of the protein (Fig.\ref{figHIU})(c)):
\begin{align*}
 U-U-U&-D-U-U -U-U-D -D-U-U -D-D-U \\
& -D-U-U -U-U-D -D-U-U -D-D-D -D-U-U \\
& -U-D-D -D-D-D -D-U-U -D-D-U -D-U-U \\
& -U-U-D -D-U-U -U-U-D -D-U-U -U-U-D \\
& -D-D-D.
\end{align*}
Table \ref{tab3} shows the correspondence between the $U/D$ code and 
the amino-acid sequence of the protein. (Also see Fig.\ref{figHIU}(c) and (d).)

\begin{table}[tb]
\caption{$U/D$ code and the amino-acid sequence of 2HIU chain A. 
($0$ denotes $D-D-D$ , $1$ denotes $D-D-U$ and so on.) }
\label{tab3}
\renewcommand\arraystretch{1.5}
\noindent\[
\begin{array}{c||cccc cccc}
{\text{No.}} &{1}&{2}&{3} &{4}&{5}&{6}&{7}&{8} \\
\hline
{\text{Amino-acid}} &{\text{GLY}}&{\text{ILE}}&{\text{VAL}} &{\text{GLU}} 
&{\text{GLN}}&{\text{CYS}}  &{\text{CYS}}&{\text{THR}}  \\
{\text{$U/D$ code}} &{7}&{3}&{6}  &{3}&{1}&{3}  &{6}&{3}  \\
\hline
\hline
{\text{No.}} &{9}&{10}&{11} &{12}&{13}&{14}&{15}&{16} \\
\hline
{\text{Amino-acid}} &{\text{SER}}&{\text{ILE}}&{\text{CYS}} &{\text{SER}}
&{\text{LEU}}  &{\text{TYR}}&{\text{GLN}}&{\text{LEU}}  \\
{\text{$U/D$ code}} &{0}&{3}&{4} &{0} &{3}&{1}&{3}&{6}  \\
\hline
\hline
{\text{No.}} &{17}&{18}&{19} &{20}&{21}&{}&{}&{} \\
\hline
{\text{Amino-acid}} &{\text{GLU}}&{\text{ASN}}&{\text{TYR}} &{\text{CYS}}
&{\text{ASN}}  &{}&{}&{}  \\
{\text{$U/D$ code}} &{3}&{6}&{3} &{6} &{0}&{}&{}&{}  \\
\hline
\end{array}
\]
\end{table}

\appendix

\renewcommand{\thesection}{\Alph{section}}
\section{Differential geometry of $N$-hedron tiles}

\subsection{Space of $N$-hedron tiles}

Let $L_N{}^\ast $ be the collection of all integer points of the $N$-dimensional Euclidean space $\mathbb{R}^N$:
\[
L_N{}^\ast :=\left\{x_1{}^{l1} x_2{}^{l2} \cdots x_N{}^{lN}\ |\  l_i \in \mathbb{Z} 
\text{ for all } i \right\}.
\]
And consider the collection $S$ of all ``slant'' $N$-hedrons defined by $L_N{}^\ast$:
\[
S := \left\{ a\left[x_{\rho(1)} \cdots x_{\rho(N-1)}\right] 
            \ |\   a \in L_N{}^\ast, \  \rho \in \mathit{S_N} \right\},
\]
where $\mathit{S_N}$ is the $N$-th symmetric group and 
$a\left[x_{\rho(1)} \cdots x_{\rho(N-1)}\right]$ denotes the convex hull 
$conv[a_0, a_1, \ldots, a_{N-1}]$ of $N$ points 
$a_0=a, a_1=ax_{\rho(1)}, \ldots, a_{N-1}=ax_{\rho(1)}x_{\rho(2)}$ 
$\cdots x_{\rho(N-1)}$ in $\mathbb{R}^N$: 
\[
a\left[x_{\rho(1)} \cdots x_{\rho(N-1)}\right]:= 
\left\{ \prod_{0 \le i <N}  a_i^{\lambda_i}\ |\  
 0\le \lambda_i \in \mathbb{R} \text{ s.t. }\sum_{0 \le i <N} \lambda_i=1   \right\}.
\]

The collection $B$ of all ``flat'' $N$-hedrons is defined as the quotient of $S$ by 
``shift operator'' $\sigma$ on $S$ (Fig.\ref{figAPP}(a)). That is, $B:= S / \sigma$, where
\[
\sigma\left(a\left[x_{\rho(1)} \cdots x_{\rho(N-1)}\right]\right) 
:= ax_{\rho(1)}\left[x_{\rho(2)} \cdots x_{\rho(N)}\right]. 
\]

\subsection{Differential structure on $B$}

``Tangent bundle'' $T[B]$ on $B$ is defined as the quotient of $S$ by $\sigma^N$:
\begin{align*}
&T[B]:= S/ \sigma^N, \\ 
&\pi: T[B] \to B, \  \pi\left(s \mod \sigma^N\right):= s \mod \sigma.
\end{align*}

We identify $T[B]$ with $B \times \{e/x_1, e/x_2, \ldots, e/x_N \}$ ($e=x_1x_2 \cdots x_N$) 
by one-to-one correspondence
\[
s \mod \sigma^N  \sim ( s \mod \sigma, Ds),
\]
where the ``gradient'' $Ds$ of $s \in S$ is defined by
\[
D a\left[x_{\rho(1)} \cdots x_{\rho(N-1)}\right] :=x_{\rho(1)} \cdots x_{\rho(N-1)}=e/ x_{\rho(N)}.
\]

Let $s=a\left[x_{\rho(1)} \cdots x_{\rho(N-1)}\right] \in S$. 
Then $s \mod \sigma^N \in T[B]$ specifies ``local trajectory'' 
$\left\{s_U \mod \sigma, s \mod \sigma, s_D \mod \sigma \right\}$ 
at $s \mod \sigma \in B$ (Fig.\ref{figAPP}(b)), where  
\begin{align*}
&s_U:=a\left[x_{\rho(1)} \cdots x_{\rho(N-2)} x_{\rho(N)}\right], \\
&s_D:=ax_{\rho(1)}\left[x_{\rho(2)} \cdots x_{\rho(N-1)} x_{\rho(1)}\right].
\end{align*}
And we shall obtain a flow on $B$ by patching these local trajectories together.

\subsection{Cones and their boundary surfaces}

Let $\mathbb{PHN}^N:=\{Cone^\ast A \ | \ A \subset L_N{}^\ast \}$, where
\[
Cone^\ast A:=\left\{p x_1{}^{l1} x_2{}^{l2} \cdots x_N{}^{lN}\in L_N{}^\ast \ | \ 
         p \in A \text{ and } 0\le l_i \in \mathbb{Z} \text{ for all } i \right\}.
\]
That is, $\mathbb{PHN}^N$ is the collection of all ``cones'' defined by $L_N{}^\ast$. 
And we denote the ``boundary surfaces'' of $w \in \mathbb{PHN}^N$ by $d_Sw$:
\[
d_S w:=\left\{ conv[a_0, a_1, \ldots, a_{N-1}] \in S \ | \  \text{$l_w(a_i) = 0$ for all $i$} \right\},
\]
where $l_w(z):= \max_{p\in w} \left\{   \min_{1\le i \le N} 
    \left\{ l_i \in \mathbb{Z} \ | \  \prod_{1\le i \le N}  y_i{}^{li}=z/p \right\} \right\}$ for $z \in L_N{}^\ast$.

The boundary surfaces of a cone induce a vector field on $B$. 

\subsection{Vector field on $B$}

Let $w \in \mathbb{PHN}^N$. Then $d_Sw$ specifies a unique $N$-hedron $s \in d_Sw$ over 
each $t \in B$, which we denote by $\Gamma_w(t)$:
\[
\Gamma_w(t):= \text{ the unique $N$-hedron $s \in d_Sw$ s.t. $t=s \mod \sigma$}.
\]
And $\Gamma_w$ induces vector field $X_w$ over $B$: 
\[
 X_w(s \mod \sigma ):=D\Gamma_w(s \mod \sigma ). 
\]

\begin{figure}[tb]
\includegraphics{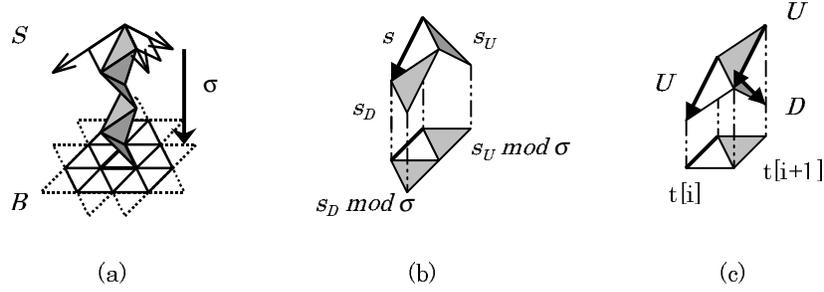}
\caption{Differential geometry of $3$-hedron tiles. 
(a): Fiber of $S$ over a point of $B$. 
(b): The local trajectory specified by $s \in S$. 
(c): The second derivative along orbit $\{ t[i] \}$.}
\label{figAPP}
\end{figure}

Let $\{ t[i] \} \subset B$ be a trajectory defined by vector field $X_w$. 
And we define the ``second derivative'' $D^2 \Gamma_w(t[i])$ of $\Gamma_w$ along 
$\{t[i]\}$ as a $\{ U, D \}$-valued function by
\[
D^2 \Gamma_w(t[i+1]):=  \begin{cases}
                   D^2 \Gamma_w(t[i])  \quad  \text{if  $X_w(t[i+1])=X_w(t[i])$}, \\
                   - D^2 \Gamma_w(t[i])  \quad  \text{else},
                   \end{cases}
\]
where $- D:= U$ and $-U:=D$ (Fig.\ref{figAPP}(c)).

Then we can encode the $N-1$-dimensional structure of any trajectory by the second 
derivative along the trajectory, i.e., an $U/D$ sequence.

\bibliographystyle{amsplain}

\end{document}